\documentclass{article}

\usepackage{amssymb}
\usepackage{amsmath,amsthm,amsfonts}

\newtheorem{theorem}{Theorem}
\newtheorem{proposition}{Proposition}
\newtheorem{lemma}{Lemma}
\newtheorem{corollP}{Corollary}[proposition]
\theoremstyle{remark}
\newtheorem{remark}{Remark}

\newtheorem{definition}{Definition}

\newcommand{\Naturales}{\mathbb{N}}
\newcommand{\F}{{\mathcal F}}
\newcommand{\e}{\epsilon}
\begin{document}

\title{Measures related to $(\epsilon,n)$-complexity functions}

\author{Valentin Afraimovich and Lev Glebsky}
\thanks{The authors were partially supported by CONACyT grant C02-42765, L.G. was partially supported by
CONACyT grant 50312}

\date{}
\maketitle
\begin{abstract}
The $(\epsilon,n)$-complexity functions describe total instability
of trajectories in dynamical systems. They reflect an ability of trajectories
going through a Borel set to diverge on the distance $\epsilon$ during
the time interval $n$. Behavior of the $(\epsilon,
n)$-complexity functions as $n\to\infty$ is reflected in the
properties of special measures. These measures are constructed as
limits of atomic measures supported at points of
$(\epsilon,n)$-separated sets. We study such measures. In
particular, we prove that they are invariant if the
$(\epsilon,n)$-complexity function grows subexponentially.

{\bf keywords} Topological entropy, complexity functions,
separability.

AMC: 28C15, 37C99

\end{abstract}

\section{Introduction}

The instability of orbits in dynamical systems is quantitatively reflected
by complexity functions. Topological complexity reflects pure topological features
of dynamics \cite{BHM}, symbolic complexity (see, for instance, \cite{Fe})
deals with symbolic systems, and the $(\epsilon,n)$-complexity (see definition below)
depends on a distance in the phase space. If a dynamical system is generated by a
map $f:X\rightarrow X$ where $X$ is a metric space with a distance $d$,
one can introduce the sequence of distances (\cite{Bo})
\begin{displaymath}
d_n(x,y)=\max_{0\leq i\leq n-1} d(f^i x,f^i y), n\in \mathbb{N},
\end{displaymath}
and study the $\epsilon$-complexity with respect to the distance $d_n$ as a function
of ``time'' $n$. This function describes the evolution of instability of orbits in time
(see, for instance, \cite{AZ},\cite{ZA},\cite{AG},\cite{AG1}). It depends not only on $n$
but on $\epsilon$ as well.\\

In fact, the $(\epsilon,n)$-complexity $C_{\e,n}$ is the maximal number of $\epsilon$-distinguishable pieces of
trajectories of temporal length $n$. It is clear that this number is growing as $\epsilon$ is
decreasing. If a system possesses an amount of instability then this number is also
growing as $n$ is increasing: the opportunity for trajectories to diverge on the distance
$\epsilon$ during the temporal interval $n+1$ is greater than to do it during $n$ temporal
units. It is known (see, for instance, \cite{T}) that
$$
b:=\lim_{n\to\infty}\limsup_{\e\to 0} \frac{\ln C_{\e,n}}{-\ln \e}
$$
is the fractal (upper box) dimension of $X$. Moreover (see, for instance, \cite{Bo1})
$$
h:=\lim_{\e\to 0}\limsup_{n\to\infty}\frac{\ln C_{\e,n}}{n}
$$
is the topological entropy of the dynamical system $(X,f)$. Following Takens \cite{T},
people say that the system $(X,f)$ is deterministic and possesses dynamical chaos if
$0<b<\infty$, $0<h<\infty$.

Thus, the $(\epsilon,n)$-complexity is a global characteristic of the evolution of
instability allowing one to single out systems with dynamical chaos.

Generally, this process of evolution occurs very non-uniformly: there are regions in the phase
space (red spots) where the divergence of trajectories is developing very fast, for
other pieces of initial conditions (green spots) trajectories manifest the distal behavior
for long intervals of time and only after that have a possibility to diverge on the
distance $\epsilon$. In other words, for every fixed $n$ there is a ``distribution'' of
initial points according to their ability to diverge to the distance $\epsilon$ during
the time interval $n$.

In our article we prove that these distributions converge to the limiting ones as
$n\to\infty$. We call them the measures related to the $(\epsilon,n)$-complexity.
We study main properties of the measures and consider some examples that allow us
to hypothesize that generally for systems with positive topological entropy these
measures are non-invariant. We also prove that for systems with zero topological
entropy these measures are invariant.

The article is a continuation of the study we started in our previous work \cite{AG} where we
proved the existence of special measures reflecting the asymptotic behavior of the
$(\epsilon,n)$-complexity as $\epsilon \rightarrow 0$.\\

\section{Set-up and definitions}

Let $(X,d)$ be a compact metric space with a distance $d$, $S\subset X$ and
$f:X\backslash S\rightarrow X$ a continuous map. If $S=\varnothing$, the map is continuous on $X$;
if not, it can be discontinuous. We assume that $X\backslash S$ is open dense in $X$ set.
Further properties of $S$ will be specified below.\\

Let $\mathcal{D}=\displaystyle{\bigcup^{\infty}_{i=0}} f^{-i}S$, the set of all preimages
of $S$. The dynamical system $(f^{i}, X\backslash \mathcal{D})$ and the distances $d_n$
are well-defined.\\

The notion of the $\epsilon$-separability was first introduced by Kolmogorov and Tikhomirov
\cite{KT} and was applied to study dynamical systems by Bowen \cite{Bo}.\\

\begin{definition}\label{Def1}
Two points $x$ and $y$ in $X\backslash \mathcal{D}$ are said to be
$(\epsilon,n)$-separated if $d_n(x,y)\geq \epsilon$.
\end{definition}

It means that the pieces of the orbits of temporal length $n$ going through $x$ and $y$
are $\epsilon$-distinguishable at the instant $i$ of time, $0\leq i\leq n-1$.\\

\begin{definition}\label{Def2}
\begin{itemize}
\item[(i)] A set $Y\subset X\backslash \mathcal{D}$ is said to be $(\epsilon,n)$-separated if
any pair $x,y\in Y$, $x\neq y$, is $(\epsilon,n)$-separated.
\item[(ii)] Given $A\subset X\backslash \mathcal{D}$, the quantity
$C_{\epsilon,n}(A)=\max \{|Y|:Y\subset A \mbox{ is }
(\epsilon,n)-\mbox{separated}\} $
where $|Y|$ is the cardinality of $Y$, is called the $(\epsilon,n)$-complexity of the
set $A$. As the function of $n$ it is called the $\epsilon$-complexity function of $A$.
\end{itemize}
\end{definition}

The following proposition is proved exactly en the same way as in \cite{AG}.\\

\begin{proposition}\label{Prop1}
Given $B_1,B_2\subset X\backslash \mathcal{D}$ and $\epsilon > 0$, the following
inequality holds
\begin{displaymath}
C_{\epsilon,n}(B_1 \cup B_2)\leq C_{\epsilon,n}(B_1) + C_{\epsilon,n}(B_2).
\end{displaymath}
\end{proposition}

\begin{definition}\label{Def3}
Given $Z\subset X\backslash \mathcal{D}$, an $(\epsilon,n)$-separated set $Y\subset Z$
is called $(\epsilon,n)$-optimal in $Z$ if $|Y|=C_{\epsilon,n}(Z)$.
\end{definition}

\section{Measures}

Here we define some measures reflecting the asymptotic behavior of the
$(\epsilon,n)$-complexity as $n\rightarrow \infty$. For that, we use a technique
of ultrafilters (see Appendix 1) and also the Marriage Lemma (Appendix 2).\\

Given $\epsilon > 0$, $n\in \mathbb{Z}$, consider an optimal (in $X\backslash \mathcal{D}$)
$(\epsilon,n)$-separated set $A_{\epsilon,n}$. Allowing $n\rightarrow \infty$ we fix a
sequence of $(\epsilon,n)$-separated sets. Introduce the following functional\\
\begin{displaymath}
I_{\epsilon,n}(\phi)= \frac{1}{C_{\epsilon,n}}\displaystyle\sum_{x\in A_{\epsilon,n}}\phi(x).
\end{displaymath}
where $\phi: X\rightarrow \mathbb{R}$ is a continuous function. It is clear that $I_{\epsilon,n}$
is a positive bounded linear functional on $C(X)$. Moreover, for any fixed $\epsilon > 0$
the sequence $I_{\epsilon,n}(\phi)$ is bounded. Fix an arbitrary non-proper ultrafilter
$\mathcal{F}$, see Appendix I. Consider
\begin{displaymath}
I_\epsilon(\phi)= \displaystyle\lim_{\mathcal{F}}I_{\epsilon,n}(\phi).
\end{displaymath}

$I_{\epsilon}$ is a positive bounded linear functional on $C(X)$ that may depend on the choice of $\epsilon$,
the ultrafilter $\mathcal{F}$ and optimal sets $A_{\epsilon,n}$. We denote by
$\mu=\mu_{\epsilon,\mathcal{F},A_{\epsilon,n}}$ the corresponding regular Borel measure on
$X$.\\
{\it Remark.} As one can see, the functional $I_\epsilon (\cdot)$ is defined for any bounded function,
in particular, for the characteristic function $\chi_Y$ of a set $Y$. Generally,
$I_\epsilon(\chi_Y)\neq\mu(Y)$. But if $C$ is a compact set and $W$ is an open set then
$I_\epsilon(\chi_C)\leq\mu(C)$ and $I_\epsilon(\chi_W)\geq\mu(W)$, see \cite{GG,Halm}.

\begin{definition}\label{Def4}
The measures $\{\mu\}$ will be called the measures related to the $(\epsilon,n)$-complexity.
\end{definition}

In the further consideration we will use the following property of
a measure $\mu$.

\begin{proposition}\label{Prop3'}
If $\mu(S)=0$ then for any sequence of positive numbers $\delta_n$,
$\delta_n\to 0$ as $n\to\infty$, one has
\begin{equation}\label{Eq1'}
\lim_{\mathcal F}\frac{1}{C_{\epsilon,n}}|A_{\epsilon,n}\setminus
O_{\delta_n}(S)|=1,
\end{equation}
where $A_{\epsilon,n}$ are the $(\epsilon,n)$-optimal sets used in
the definition of $\mu$ and $O_{\delta}(S)$ is the
$\delta$-neighborhood (in the metric $d$) of the set $S$.
\end{proposition}
\begin{proof}
In fact, the validity of (\ref{Eq1'}) follows directly from the
definition of $\mu$. Indeed, for any small $\delta>\alpha>0$
$$
\lim_{\mathcal F}\frac{1}{C_{\epsilon,n}}|A_{\epsilon,n}\cap
O_{\alpha}(S)|=I_\epsilon(\chi_{O_\alpha(S)})\leq
                                       I_\epsilon(\chi_{\overline{O_\alpha(S)}}) \leq\mu(\overline{O_\alpha(S)})
\leq\mu(O_\delta(S)),
$$
and $\mu(O_\delta(S))\to 0$ as $\delta\to 0$ ($\mu$ is a regular
measure). Moreover, $\delta_n<\delta$ if $n$ is large enough.
Therefore,
$$
\lim_{\mathcal F}{1\over C_{\epsilon,n}}|A_{\epsilon,n}\cap
O_{\delta_n}(S)|\leq\mu(O_\delta(S)).
$$
It implies the desired result.
\end{proof}
The following proposition is proved in the same way as Proposition 6 in \cite{AG},
(one should just replace the distance $d$ by the distance $d_n$ and apply the
Marriage Lemma (see Appendix II). For completeness we present the proof here. \\

\begin{proposition}\label{Prop3}
Let $A$ and $B$ be the $(\epsilon,n)$-separated sets and $A$ is
optimal. Then there exists an injection map $\alpha_n:B \rightarrow
A$ such that $d_n(x,\alpha_n(x))< \epsilon$ for any $x\in B$. If
$|B|=|A|$ then $\alpha_n$ is bijection.
\end{proposition}
\begin{proof}
Recall that $O_\epsilon (x)=\{y\; :\; d_n(x,y)<\epsilon\}$, the ball of radius $\epsilon$
centered at $x$. Given $Y\subseteq X$ let $O_\epsilon(Y)=\bigcup\limits_{x\in Y} O_\epsilon(x)$.

For any $x\in B$ let $A_x=O_\epsilon (x)\cap A$. If we show that for any $S\subseteq B$ the
following inequality holds
\begin{equation} \label{Hall_condition}
|\bigcup_{x\in S} A_x|\geq |S|,
\end{equation}
then the proposition follows from the Marriage lemma.
To prove the inequalities~(\ref{Hall_condition}), suppose that
$|\bigcup\limits_{x\in S} A_x|=|O_\epsilon(S)\cap A|< |S|$ for some $S\subseteq B$. Then
$$
|S\cup(A\setminus (O_\epsilon(S)\cap A)|=|S|+(|A|-|O_\epsilon(S)\cap A|)>|A|=C_{\e,n}.
$$
On the other hand, the set $S\cup(A\setminus (O_\epsilon(S)\cap A)$ is $(\e,n)$-separated.
We have a contradiction with optimality of $A$.
\end{proof}

For an arbitrary map $f$ the validity of the inequalities $d_n(x,\alpha_n(x))< \epsilon$,
$n\in \mathbb{N}$, does not imply that $d(x,\alpha_n(x))\rightarrow 0$. For example,
for distal dynamical system it is not true. As a corollary, we have an unpleasant fact
that the functional $I_{\epsilon}$ and the corresponding measure $\mu$ may depend on the
choice of optimal sets. In the next section we introduce a class of maps for which it is
not so.

\section{Measures for $\epsilon$-expansive maps}

We begin with the following definition.\\

\begin{definition}\label{Def5}
\begin{itemize}
\item[(i)] We say that the map $f$ is $\epsilon$-expansive if for any $\delta > 0$ and
any pair $x,y\in X\backslash \mathcal{D}$, $x\neq y$, there exists $N=N(x,y,\delta)$
such that the inequality $d_n(x,y)\leq \epsilon$, $n\geq N$, implies that
$d(x,y)\leq \delta$.
\item[(ii)] The map $f$ is uniformly $\epsilon$-expansive if there exists a sequence of
non-negative numbers $\delta_n \rightarrow 0$ as $n\rightarrow
\infty$ such that for any pair $x,y\in X\backslash \mathcal{D}$ with
$d_n(x,y)\leq \epsilon$ one has $d(x,y)\leq \delta_n$, $n=1,2,\dots$
\end{itemize}
\end{definition}

\begin{lemma}\label{Lem1}
A continuous $\epsilon$-expansive map ($S=\varnothing$) is uniformly $\epsilon$-expansive.
\end{lemma}

\begin{proof}
Assume that it is not true, i.e. there exists a sequence $n_k \rightarrow \infty$ as
$k \rightarrow \infty$ and a sequence of pairs $x_k \neq y_k$ such that
$d_{n_k}(x_k,y_k)\leq \epsilon$, $d(x_k,y_k) > \beta > 0$. Since $X$ is compact, then
without loss of generality one may assume that there exist
$\displaystyle\lim_{k\rightarrow \infty}x_k=x_0$,
$\displaystyle\lim_{k \rightarrow \infty}y_k=y_0$ (in the metric $d$) and $d(x_0,y_0)\geq \beta$.
Since $d_{n_k}(x_k,y_k)\leq \epsilon$, then $d(x_k,y_k)\leq \epsilon$ and $d(x_0,y_0)\leq \epsilon$.
Also, $d(fx_k,fy_k)\leq \epsilon$, so, $d(fx_0,fy_0)\leq \epsilon$, because of the continuity of
$f$. In the same way, one may show that $d(f^{m-1}x_0,f^{m-1}y_0)\leq \epsilon$
(if one chooses $n_k > m$), thus $d_m(x_0,y_0)\leq \epsilon$ for any $m\in \mathbb{N}$.
Since $f$ is $\epsilon$-expansive $x_0 = y_0$, a contradiction.

\end{proof}

For uniformly $\epsilon$-expansive maps the following fact takes place.\\

\begin{proposition}\label{Prop4}
If $f$ is uniformly $\epsilon$-expansive, $y,z\in f^{-1}x$, $x\in X\backslash \mathcal{D}$, and
$d(z,y) \leq \epsilon$ then $z=y$.
\end{proposition}

\begin{proof}
Since $fz = fy$ then $d_k(z,y)\leq \epsilon$ for every $k \in \mathbb{N}$.
Hence, $d(z,y)\leq \delta_k$, i.e. $z=y$.

\end{proof}

\begin{theorem}\label{Th1}
If $f$ is uniformly $\epsilon$-expansive then the functional $I_\epsilon$
(and the corresponding measure) is independent of the choice of optimal sets $A_{\epsilon,n}$.
\end{theorem}

\begin{proof}
Let $A_{\epsilon,n}$, $B_{\epsilon,n}$ be optimal $(\epsilon,n)$-separated sets,
$n\in \mathbb{N}$. Because of Proposition~\ref{Prop3}, there exists a bijection
$\alpha_n:A_{\epsilon,n} \rightarrow B_{\epsilon,n}$ such that
$d_n(x,\alpha_n(x))< \epsilon$. It implies the existence of $\delta_n\geq 0$
such that $d(x,\alpha_n(x)) \leq \delta_n$. Thus,
\begin{displaymath}
\frac{1}{C_{\epsilon,n}} \Big| \displaystyle\sum_{x\in A_{\epsilon,n}}\phi(x) - \displaystyle\sum_{x\in B_{\epsilon,n}}\phi(x)\Big| \leq \frac{1}{C_{\epsilon,n}} \displaystyle\sum_{x\in A_{\epsilon,n}}
\big|\phi(x) - \phi(\alpha_n(x))\big| \leq \omega_{\delta_n}(\phi),
\end{displaymath}
where
$$
\omega_{\delta_n}(\varphi)=\sup_{d(x,y) \leq \delta_n}\mid
\varphi(x)-\varphi(y)\mid,
$$
the modulus of continuity of $\phi$.\\

Since $X$ is compact and $\phi$ is continuous, then $\omega_{\delta_n}(\phi)\rightarrow 0$
as $\delta_n \rightarrow 0$.

\end{proof}

\section{Non-invariance of the measures}
 As it was mentioned in Introduction, we have studied in
 \cite{AG}
 behavior of $C_{\e,n}$ as $\e\to 0$. In particular, we proved that
 for any sequence $\e_k\to 0$, $k\to\infty$, there exists a regular Borel
 measure corresponding to the functional $I(\cdot)=\lim\limits_\F
 I_{\e_n,1}$, where $\F$ is a nonproper ultrafilter. We call here
 such measures the $\e$-measures. The measures constructed in
 Section~4 will be called the $n$-measures.
Neither $\e$-measures nor $n$-measures  are not obliged to be
$f$-invariant. The following example shows that it is really so. In
the example $n$-measure will coincide with $\epsilon$-measure.

The item (ii) of Definition~\ref{Def5} could be rewritten as
"$d(x,y)>\delta_n$ implies $d_n(x,y)>\epsilon$. It means that any
$(\delta_n,1)$-separated set is an $(\e,n)$-separated set. So,
suppose that $f$ is continuous and satisfies the following stronger
condition:
\begin{itemize}
\item[(ii*)] There exists a sequence of
non-negative numbers $\delta_n \rightarrow 0$ as $n\rightarrow
\infty$ such that for any pair $x,y\in X$ one has
$d_n(x,y)\leq \epsilon$ if and only if $d(x,y)\leq \delta_n$,
$n=1,2,\dots$
\end{itemize}
In this case set $A$ is $(\e,n)$-separated if and only if it is
$(\delta_n,1)$-separated and, $A$ is an optimal $(\e,n)$-separated
set if and only if it is an optimal $(\delta_n,1)$-separated set.
So, the $\e$-measure, corresponding to $\lim\limits_\F
 I_{\delta_n,1}$   equals the $n$-measure, corresponding to $\lim\limits_\F
 I_{\e,n}$.

It easy to check that any symbolic dynamical system with finite alphabet
satisfies (ii*), but
the corresponding $\epsilon$-measure, constructed in \cite{AG} is
not shift-invariant. Let us describe an example.\\
{\it Example.} Let $X=\Omega_{M}$ be a topological Markov chain,
defined by a finite matrix $M:\{0,1,...,p-1\}^2\to\{0,1\}$, i.e.
$\Omega_{M}=\{(x_0,x_1,...)\;|\; x_i\in
\{1,2,...,p-1\}\;\mbox{and}\; M(x_i,x_{i+1})=1\}$. The set $\Omega_M$ is endowed
with the  metric
$d(x,y)=2^{-n}$, where $n=\min\{i\;|\;x_i\neq y_i\}$, and map $f$ is
the shift: $f(x_0,x_1,x_2,...)=(x_1,x_2,...)$. The map $f$ satisfies (ii*)
with $\delta_n=2^{-n+1}\e$. So, here $n$-measure and $\e$-measures
coincide and are given by the following proposition (see
\cite{AG}).
\begin{proposition}
Let $M$ be a primitive matrix and $C\subset \Omega_M$ is an
admissible cylinder of length $n$, ending by $i$. Then
$\mu(C)=\lambda^{-n}e_i$, where $(e_0,e_1,...,e_{p-1})$ is
the positive eigenvector of $M$, with $e_0+e_1+...e_{p-1}=1$,
and $\lambda>1$ is the corresponding eigenvalue.
\end{proposition}
If $\lambda$ is not an integer such a measure can not be
invariant.

  Other properties of the measures we want to discuss here are related to
the group $GI(X)$ of isometries of $X$.

\begin{definition}
We will say that $f$ commute with a group $G$ of transformations of
$X$ iff for any $g\in G$ there exists $h\in G$ such that $f\circ g
=h\circ f$.
\end{definition}
{\it Example.} The map $f:x\to 2x, \mod{1}$, of the circle commutes with the group of rotations
$x\to x+\omega, \mod{1}$.

One can check that if $f$ commutes with the group $GI(X)$ of
isometries of $X$, then elements of $GI(X)$ are isometries for
$d_n(\cdot,\cdot)$, and the following proposition is true:
\begin{proposition}
If $f$ commutes with $GI(X)$ then the  corresponding $n$-measure is
$GI(X)$ invariant.
\end{proposition}
It follows that in the example above the $n$-measure is just the Lebesgue measure.

Let us present an easy example where the $\e$-measure is different
from $n$-measure. The example in some sense is artificial but it
shows that if the  expansivity of a  map is different at different points, then the $n$-measure may be
different from the $\e$-measure. Let $X=\{0,1\}\times [0,1)$ (disjoint
union of to unit intervals, considered as a circles,
$d((i,x),(i,y))=\min\{|x-y|,||x-y|-1/2|\}$, $d((0,x),(1,y))=1$.
Let
$f(0,x)=(0,\;2x\mod 1)$, $f(1,x)=(1,\;3x\mod 1)$. Then the $\e$-measure
is a Lebesgue measure, such that $\mu(\{0\}\times
[0,1))=\mu(\{1\}\times [0,1))=1/2$, for the $\e$-measure of the circle is a Lebesgue measure and
$C_\e(\{0\}\times [0,1))=C_\e(\{1\}\times [0,1))$  But the $n$-measure
$\mu(\{0\}\times [0,1))=0$, since $C_{\e,n}(\{0\}\times [0,1))/C_{\e,n}(\{1\}\times [0,1))\to 0$ as
$n\to\infty$.

 We think that generally for dynamical systems with positive
topological entropy, $n$-measures are not invariant. But for systems
with zero entropy, they may be invariant in a general enough
situation.

\section{Invariance of the measures}

For many subexponential functions $C_{\epsilon,n}$, the following equality holds

\begin{equation}\label{Eq1}
\displaystyle\lim_{n\rightarrow \infty}\frac{C_{\epsilon,n}-C_{\epsilon,n-1}}{C_{\epsilon,n}}=0.
\end{equation}

\begin{remark}
It follows from the definition of the topological entropy that the equality (\ref{Eq1}) is not satisfied if
the topological entropy $h_{top}(f|X\backslash\mathcal{D}) > 0$ and $\epsilon$ is small enough.
\end{remark}

In fact, (\ref{Eq1}) could not be satisfied even if
$h_{top}=0$. Suppose that $C_{\epsilon,n}=2^{\big[\sqrt{n}\big]}$,
where $\big[{\mathbf\cdot}\big]$ means the integer part of the
number. For this sequence
$$
\frac{C_{\epsilon,n}-C_{\epsilon,n-1}}{C_{\epsilon,n}}=\left\{\begin{array}{cc}
                                      \frac{1}{2}, & \mbox{if } n=m^2 \\
                                       0, & \mbox{if $n$ is not a full square}
                                       \end{array}\right.
$$
So,  for subexponential functions $C_{\epsilon,n}$ limit (\ref{Eq1})
could not exist. But, for any subexponential $C_{\epsilon,n}$ there  exists
the lower limit:
$$
\displaystyle\liminf_{n\rightarrow
\infty}\frac{C_{\epsilon,n}-C_{\epsilon,n-1}}{C_{\epsilon,n}}=0,
$$

(it equals $0$ since if the lower limit $>0$ then $C_{\epsilon,n}$ grows exponentially,
the contradiction).
It implies that there exists an ultrafilter such that the
corresponding limit with respect to this ultrafilter is $0$. So, we
replace (\ref{Eq1}) by the following more weak assumption:
\begin{equation}\label{Eq2}
\displaystyle\lim_{\mathcal{F}} \frac{1}{C_{\epsilon,n}} \big(C_{\epsilon,n}-C_{\epsilon,n-1}\big)=0,
\end{equation}
where $\mathcal{F}$ is a non-proper ultrafilter.\\

The assumption (\ref{Eq2}) imply

\begin{proposition}\label{Prop5}
Let
\begin{displaymath}
\{n\;| \;C_{\epsilon,n-1} \leq b_n \leq
C_{\epsilon,n},\}\in\mathcal{F}.
\end{displaymath}
Then
\begin{equation}\label{Eq3}
b_n = C_{\epsilon,n}(1-q_n)
\end{equation}
where $\displaystyle\lim_\mathcal{F}q_n=0$.
\end{proposition}

\begin{proof}
Defining $q_n$ as $q_n:= 1- \frac{b_n}{C_{\epsilon,n}}$ we need to show only that
$\displaystyle\lim_{\mathcal{F}} q_n= 0$. Assume not, i.e.
$\displaystyle\lim_{\mathcal{F}} q_n =\rho > 0$. Then
\begin{displaymath}
C_{\epsilon,n}-b_n = q_n C_{\epsilon,n} = C_{\epsilon,n}(\rho + \xi_n)
\end{displaymath}
with $\displaystyle\lim_{\mathcal{F}} \xi_n= 0$. Thus,
\begin{displaymath}
b_n= C_{\epsilon,n}-q_nC_{\epsilon,n} = C_{\epsilon,n} (1-\rho-\xi_n) \geq C_{\epsilon,n-1}.
\end{displaymath}
Hence,
\begin{equation}\label{Eq4}
\frac{C_{\epsilon,n} - C_{\epsilon,n-1}}{C_{\epsilon,n}} \geq \rho + \xi_n,
\end{equation}
a contradiction with (\ref{Eq2}).

\end{proof}

\begin{corollP}
In particular if
\begin{equation}\label{Eq5}
\{n\; |\;C_{\epsilon,n-1}=C_{\epsilon,n}(1-q_n)\}\in \mathcal{F}.
\end{equation}
then $\displaystyle\lim_\mathcal{F}q_n=0$.
\end{corollP}

>From now on we assume that $f$ is uniformly $\epsilon$-expansive.\\

\begin{proposition}\label{Prop6}
If $A_{\epsilon,n-1}$ is $(\epsilon,n-1)$-separated then $f^{-1}A_{\epsilon,n-1}$
is $(\epsilon,n)$-separated, $n=2,3,\dots$
\end{proposition}

\begin{proof}
If $y,z \in f^{-1}x$, $x\in A_{\epsilon,n-1}$, $y\neq z$, then, because of
Proposition \ref{Prop4}, $d(y,z)> \epsilon$, and $y$ and $z$ are
$(\epsilon,n)$-separated.\\

If $y$ and $z$ belong to $f^{-1}A_{\epsilon,n-1}$ and $fy \neq fz$ then they are
$(\epsilon,n)$-separated, since $fy$ and $fz$ are $(\epsilon,n-1)$-separated.

\end{proof}

\begin{corollP}
Proposition implies that $|f^{-1}A_{\epsilon,n-1}| \leq C_{\epsilon,n}$.
\end{corollP}

Let us repeat that since $f$ is uniformly $\epsilon$-expansive, the sequence $\delta_n$
is defined.\\

We restrict our attention now to a class of maps that could be
discontinuous but possess a large amount of continuity.
\begin{definition}\label{Def6} We say that $f:X \setminus S \rightarrow X$ is almost uniformly continuous if
there exist $\delta_0 >0$ such that for every $0 < \delta <
\delta_0$ and $0 < \sigma < \delta$ the
  modulus of continuity
  $$
  \omega_\sigma (f \mid X \setminus O_\delta(S))\leq \eta(\sigma)
  $$
  where the function $ \eta(\sigma)$ is independent of $\delta$ and goes to
  $0$ as $\sigma \rightarrow 0$.
\end{definition}

In other words
 $$d(fx,fy) \leq  \eta(\sigma)$$
if  $d(x,S) \geq \delta$, $d(y,S) \geq \delta$ and $d(x,y)  \leq \sigma$.

As an example, one may consider a smooth map $f$ on a subset $X \subset \textrm{R}^n$
for which $ \sup_{x \in X \setminus D} \Arrowvert \mathcal{D}f(x) \Arrowvert <
\infty$.\\

The main result of this section is the following theorem.

\begin{theorem}\label{Th2}
Let $f:X\setminus S \rightarrow X$ be an  almost uniformly
continuous, uniformly $\epsilon$-expansive map and $\mu$ be the
measure related to the $(\epsilon,n)$-complexity corresponding to
the ultrafilter $\mathcal{F}$  satisfying the equation (\ref{Eq2}).
Suppose, that $\mu(S)=0$ then $\mu$ is $f$-invariant.
\end{theorem}
\begin{proof}
If is enough to show that $I_\epsilon(\varphi)=I_\epsilon(\varphi \circ f)$
for every $\varphi \in C(X)$ where
$$
I_\epsilon(\varphi \circ f)=\lim_\mathcal{F}
\frac{1}{C_{\epsilon,n}}\sum_{x \in A_{\epsilon,n}}\varphi (fx), \hspace{2cm} A_{\epsilon, n} \subset X \setminus D.
 $$

Given an $(\epsilon,n)$-optimal $A_{\epsilon, n}$, let $A_{\epsilon, n-1}$ be an
arbitrary  $(\epsilon,n-1)$-optimal set. Then

\begin{eqnarray*}
\frac{1}{C_{\epsilon,n}} \Bigl\arrowvert \sum_{x \in A_{\epsilon,n}}\varphi (x)  &-&
\sum_{x
  \in A_{\epsilon,n}}\varphi (fx) \Bigl\arrowvert \leq          \\
\frac{1}{C_{\epsilon,n}}
  \biggl\{ \Bigl\arrowvert \sum_{x \in A_{\epsilon,n}}\varphi
  (x)-\sum_{x \in A_{\epsilon,n-1}}\varphi (x) \Bigl\arrowvert
                                                                                 &+&
\Bigl\arrowvert \sum_{x \in A_{\epsilon,n-1}}\varphi (x)-\sum_{x
  \in A_{\epsilon,n}}\varphi (fx) \Bigl\arrowvert \biggr\}
\end{eqnarray*}

{\bf The first sum}. The set $A_{\epsilon,n-1}$ is $(\epsilon, n-1)$-separated,
therefore it is $(\epsilon,n)$-separated. Proposition 3 implies that there
exists an injection
$\alpha_n: A_{\epsilon,n-1} \rightarrow A_{\epsilon,n}$ such that
$d_n(x,\alpha_n(x)) \leq \epsilon$,  and because of the uniform
$\epsilon$-expansioness of $f$, $d(x,\alpha_n(x)) \leq \delta_n$ for any $x \in
A_{\epsilon, n-1}$. Thus,

\begin{eqnarray}
\Delta^{(1)}:=\frac{1}{C_{\epsilon,n}}  \Bigl\arrowvert \sum_{x \in A_{\epsilon,n-1}}\varphi (x) &-&
\sum_{x \in A_{\epsilon,n}}\varphi (x)  \Bigr\arrowvert
\leq \nonumber \\
\frac{1}{C_{\epsilon,n}} \biggl\{ \sum_{x \in A_{\epsilon,n-1}}
  \Bigl\arrowvert \varphi
  (x)-\varphi(\alpha_n(x)) \Bigr\arrowvert
                                                                                                 &+&
 \sum_{x \in A_{\epsilon,n} \setminus
  \alpha_n(A_{\epsilon,n-1})}  \Bigl\arrowvert \varphi  (x) \Bigr\arrowvert
\biggr\}
\leq \\
 \frac{C_{\epsilon,n-1}}{C_{\epsilon,n}}
  \omega_{\delta_n}(\varphi)                                                                      &+&
\frac{C_{\epsilon,n}-C_{\epsilon,n-1}}{C_{\epsilon,n}} \cdot
 \Arrowvert \varphi \Arrowvert    \nonumber
\end{eqnarray}

where $\omega_{\delta_n}(\varphi)$ is  the modulus of continuity of $\varphi$. Since
$\varphi$ is continuous, $\omega_{\delta_n}(\varphi) \rightarrow 0$ as $\delta_n
\rightarrow 0$.\\
{\bf The second sum}. We use the identity
$A_{\epsilon,n-1}=f(f^{-1}A_{\epsilon,n-1})$ and the following representation
$$
f^{-1}A_{\epsilon,n-1}=A^{(1)}_{\epsilon,n-1} \cup A^{(2)}_{\epsilon,n-1}
$$
where $f(A^{(1)}_{\epsilon,n-1})=A_{\epsilon,n-1}$, $\mid
A^{(1)}_{\epsilon,n-1} \mid = \mid A_{\epsilon,n-1} \mid = C_{\epsilon,n-1}$
and  $A^{(2)}_{\epsilon,n-1}=f^{-1}A_{\epsilon,n-1} \setminus A^{(1)}_{\epsilon,n-1}$,
so
$\mid A^{(2)}_{\epsilon,n-1} \mid = \mid f^{-1}A_{\epsilon,n-1} \mid -
C_{\epsilon,n-1}$. Because of Proposition~\ref{Prop6},  $f^{-1}A_{\epsilon,n-1}$ is
$(\epsilon, n)$-separated, therefore, $\mid  f^{-1}A_{\epsilon,n-1} \mid \leq C_{\epsilon,n} $
and  $\mid A^{(2)}_{\epsilon,n-1} \mid \leq  C_{\epsilon,n} -
C_{\epsilon,n-1}$. Now,
\begin{equation*}
\sum_{x \in A_{\epsilon,n-1}}\varphi (x)=\sum_{x \in A^{(1)}_{\epsilon,n-1}}\varphi (fx)=
\sum_{x \in f^{-1}A_{\epsilon,n-1}}\varphi (fx)-
\sum_{x \in A^{(2)}_{\epsilon,n-1}}\varphi (fx),
\end{equation*}
hence,

 \begin{eqnarray*}
\Delta^{(2)}:=\frac{1}{C_{\epsilon,n}}  \Bigl\arrowvert \sum_{x \in A_{\epsilon,n-1}}\varphi (x)
-\sum_{x\in A_{\epsilon,n}}\varphi (fx)  \Bigl\arrowvert \leq \\  \nonumber
\frac{1}{C_{\epsilon,n}}
\left\{\Bigl\arrowvert \sum_{x \in f^{-1}A_{\epsilon,n-1}} \varphi
  (fx) - \sum_{x \in A_{\epsilon,n}} \varphi (fx)  \Bigl\arrowvert
+
  \Bigl\arrowvert  \sum_{x \in A^{(2)}_{\epsilon,n-1}}\varphi
  (fx)  \Bigl\arrowvert\right\}.  \nonumber
\end{eqnarray*}

Since $f^{-1}A_{\epsilon,n-1}$ is $(\epsilon,n)$-separated, there exists an
injection $\beta_n:f^{-1}A_{\epsilon,n-1} \rightarrow A_{\epsilon,n}$ such
that $d_n(\beta_n(x),x) \leq \epsilon$, i.e. $d(x, \beta_n(x)) \leq \delta_n$
for any $x \in f^{-1}A_{\epsilon,n-1}$. Therefore,
 \begin{eqnarray*}
\Delta^{(2)} \leq \frac{1}{C_{\epsilon,n}} \biggl\{ \sum_{x \in f^{-1}A_{\epsilon,n-1}}
   \Bigl\arrowvert \varphi (fx)-\varphi (f\beta_n(x))  \Bigr\arrowvert                        &+&
\Bigl\arrowvert \sum_{x
  \in A_{\epsilon,n} \setminus \beta_n(f^{-1} A_{\epsilon,n-1}) }\varphi (fx)  \Bigr\arrowvert \\
                                                                                              &+&
\Bigl\arrowvert \sum_{x \in A^{(2)}_{\epsilon,n-1}} \varphi (fx) \Bigr\arrowvert \biggr\}.
\end{eqnarray*}
Since $\mid  A_{\epsilon,n} \setminus \beta_n(f^{-1}  A_{\epsilon,n-1}) \mid \leq  C_{\epsilon,n} -
C_{\epsilon,n-1}$, we obtain

 \begin{equation*}
\Delta^{(2)} \leq \frac{1}{C_{\epsilon,n}} \sum_{x \in f^{-1}A_{\epsilon,n-1}}
  \Bigl\arrowvert \varphi (fx)-\varphi (f\beta_n(x))  \Bigr\arrowvert  + 2
  \frac{C_{\epsilon,n}-C_{\epsilon,n-1}}{C_{\epsilon,n}} \cdot
  \Arrowvert \varphi \Arrowvert
\end{equation*}

Because of the almost uniform continuity of $f$, we know that if $x,y \in X \setminus
O_{2\delta_n}(S)$ and $d(x,y) \leq \delta_n$, then $d(fx,fy) \leq
\eta(\delta_n)$.
Therefore
 \begin{eqnarray*}
\sum_{x \in f^{-1}A_{\epsilon,n-1}}  \Bigl\arrowvert
  \varphi(fx)                                             &-&
\varphi(f\beta_n(x))  \Bigr\arrowvert
  \leq \\
\sum_{x \in f^{-1}A_{\epsilon,n-1} \setminus O_{2 \delta_n}(S)}  \Bigl\arrowvert
  \varphi (fx)-\varphi (f\beta_n(x))  \Bigr\arrowvert
                                                                &+&
2  \Bigl\arrowvert f^{-1} A_{\epsilon,n-1}
  \cap O_{2\delta_n}(S)  \Bigr\arrowvert \cdot \Arrowvert \varphi \Arrowvert
\leq\\
\eta(\delta_n)  \Bigl\arrowvert  f^{-1}A_{\epsilon,n-1} \setminus O_{2 \delta_n}(S)
  \Bigr\arrowvert                                               &+&
2  \Bigl\arrowvert  A_{\epsilon,n} \cap O_{2\delta_n}(S)
  \Bigr\arrowvert \cdot \Arrowvert \varphi
 \Arrowvert
\leq \\
C_{\epsilon,n} \eta(\delta_n) + 2  \Bigl\arrowvert  A_{\epsilon,n}  &\cap&
 O_{2\delta_n}(S)  \Bigl\arrowvert \cdot \Arrowvert \varphi \Arrowvert
\end{eqnarray*}

and finally, we obtain

\begin{displaymath}
\Delta^{(2)} \leq \eta (\delta_n) + 2 \Arrowvert \varphi \Arrowvert \cdot \frac{\mid
  A_{\epsilon,n} \cap O_{2\delta_n}(S)\mid} {C_{\epsilon,n}} + 2 \Arrowvert
  \varphi \Arrowvert \cdot \frac{C_{\epsilon,n} -  C_{\epsilon,n-1}}{C_{\epsilon,n}}.
\end{displaymath}
Thus,
\begin{displaymath}
\Delta^{(1)} + \Delta^{(2)} \leq \omega_{\delta_n}(\varphi)+\eta (\delta_n) + 2 \Arrowvert \varphi \Arrowvert \frac{\mid
  A_{\epsilon,n} \cap O_{2\delta_n}(S)\mid} {C_{\epsilon,n}} + 3 \Arrowvert \varphi \Arrowvert \frac{C_{\epsilon,n} -  C_{\epsilon,n-1}}{C_{\epsilon,n}}
\end{displaymath}
and
$$
\lim_\mathcal{F}\biggl(\Delta^{(1)}+ \Delta^{(2)}\biggr)=0.
$$
(the third term goes to $0$ because of Proposition~\ref{Prop3'})
\end{proof}

\section{Interval exchange transformation}

Dynamical systems generated by interval exchange transformations are basic ones among
systems with zero topological entropy. They possess some amount of instability
(generally, they are weak mixing \cite{V}) and it is not difficult to calculate
their $(\epsilon,n)$-complexity functions (see below).\\

An interval exchange transformation on the interval $I=[0,1]$ can be written as
follows
\begin{displaymath}
\overline{x}=f(x) \equiv x + c_i \hspace{0.5cm} for \hspace{0.2cm} x \in [a_i,a_{i+1}),
\end{displaymath}
$i= 0,1, ..., m-1$, $a_0 = 0$, $a_m = 1$, where $c_i \neq c_{i+1}$.\\

The set of discontinuity here $S=\{ x=a_i, i=1, \dots , m-1\}$. We
assume that: $(i)$ the map $f$ is one-to-one; $(ii)$ the set
$\mathcal{D}=\displaystyle\bigcup^\infty_{k=0}f^{-k}S$ is dense in
I; $(iii)$ the set $\mathcal{D}$ does not contain $f$-periodic
points. These assumptions imply (\cite{B}, \cite{V}) that the only
invariant
measure is the Lebesgue measure.\\

\begin{proposition}
Under the assumptions $(i)$-$(iii)$, there exists $\epsilon_0 > 0$ such that
for every $0 < \epsilon < \epsilon_0$ the map $f$ is uniformly
$\epsilon$-expansive.
\end{proposition}

\begin{proof}
First of all one can find $\epsilon_0>0$ satisfying the following condition:
\begin{itemize}
\item[C1].
If $d(x,y)\leq\epsilon_0$ and interval $[x,y]$ contains a point of
discontinuity of $f$, then $d(f(x),f(y))>\epsilon_0$.
\end{itemize}
Now take $0<\epsilon\leq\epsilon_0$. It is clear that this
$\epsilon$ also satisfies the condition C1. Let
$\mathcal{D}_n=\displaystyle\bigcup^{n-2}_{k=0}f^{-k}S$, for
$n\geq2$. One can order the set
$\mathcal{D}_n=\{r_1,r_2,...,r_{k_n}\}$, $r_1<r_2<...<r_{k_n}$. Take
$\delta_1=\epsilon$ and $\delta_n=2\min\{r_{i+1}-r_i\}$. It is clear
that $\delta_n\to 0$ when $n\to\infty$.  We have to show only that
if $d(x,y)>\delta_n$ then $d_n(x,y)>\epsilon$, or
$d(f^i(x),f^i(y))>\epsilon$ for some $0\leq i <n$. But
$[x,y]\cap\mathcal{D}_n\neq\emptyset$. So, there exists $i$, $0\leq
i<n-1$ such that $[f^i(x),f^i(y)]$ contains a point of discontinuity
of $f$. But then, if $d(f^i(x),f^i(y))\leq \epsilon$ then
$d(f^{i+1}(x),f^{i+1}(y))>\epsilon$ by property C1.
\end{proof}

Thus, measure $\mu$ related to the $(\epsilon,n)$-complexity is independent of the
choice of the $(\epsilon,n)$-optimal sets. Moreover, if $\mu(S)=0$ then it is
invariant and, hence, the Lebesgue measure.\\

In fact, we can calculate the $(\epsilon,n)$-complexity function.\\

\begin{proposition}
There exists $\epsilon_0 > 0$ such that for every $0 < \epsilon < \epsilon_0$
there is $n_0 = n_0(\epsilon)$ such that for $n > n_0$.
\begin{displaymath}
C_{\epsilon,n}=(m-1)(n-1)+1.
\end{displaymath}
\end{proposition}

\begin{proof}
Let $\epsilon_0 =  \min\{\displaystyle\min_i(a_{i+1}-a_i),\displaystyle\min_j(c_{j+1}-c_{j})\}$.
This number is positive. Given $0< \epsilon < \epsilon_0$, consider the set
$A_n=\displaystyle\bigcup^{n-2}_{k=0}f^{-k}S$. Under our assumptions $(i)$ and $(ii)$,
$A_n$ does not contain the end points \{0\} and $\{1\}$ and $N_n:=|A_n|=(m-1)(n-1)$.
We denote by $b_j$, $j=1, \dots , N_n$, the points in $A_n$ ordered in such a way that
$b_0 := 0< b_1 < b_2 < \dots < b_{N_n} < 1 = : b_{N_n+1}$. The intervals
$I_j=\{b_j \leq x \leq b_{j+1}\}$, $j=0, \dots, N_n$ form a partition of $I$.
Under the assumption $(ii)$, there exists $n_0 = n_0 (\epsilon)$ such that for any
$j$, $b_{j+1} - b_j < \epsilon$ for all $n \geq n_0$.\\

Given $n\geq n_0$, consider an $(\epsilon, n)$-optimal set
$A_{\epsilon,n}= \{p_0 < p_1 < \dots < p_r\}$, $r= C_{\epsilon,n}-1$. If the pair
$(p_i,p_{i+1})$ belong to the same interval $I_j$, then
$d(f^k p_i, f^k p_{i+1})=d(p_i,p_{i+1})$ for $k=0, 1, \dots, n-1$, and we have a contradiction.
On the other side, if they belong to different intervals then there is $s$, $0 \leq s \leq n-1$,
such that $d(f^s p_i, f^s p_{i+1})\geq \epsilon_0 > \epsilon$. Thus, $C_{\epsilon,n}$ is equal
to the number of different intervals $I_j$, i.e.
\begin{displaymath}
C_{\epsilon,n}=(m-1)(n-1)+1.
\end{displaymath}

\end{proof}

\subsection{Appendix I}
Now we give some known results and definitions that can be found, for instance, in \cite{bu}.
\begin{definition}
A set $\F\subset 2^\Naturales$ is said to be a filter over $\Naturales$ iff it satisfies the following conditions:
\begin{itemize}
\item If $A\in \F$ and $B\in \F$, then $A\cap B\in\F$,
\item If $A\in \F$ and $A\subset B$ then $B\in \F$,
\item $\emptyset\not\in\F$.
\end{itemize}
\end{definition}
Let $a_n$ be a sequences of real numbers, $a$ is called to be a limit of $a_n$ with respect to
a filter $\F$, $a=\lim_\F a_n$, if for any $\e>0$ one has $\{n\ | \ |a_n-a|<\e\}\in\F$.
>From the definition of a filter it follows that $\lim_\F a_n$ is unique, if exists.\\
{\bf Example} Let $\F_F=\{A\subseteq \Naturales\ | \ \Naturales\backslash A$ is finite $\}$. $\F_F$ is said to be a Frech\'et
filter. One can check that it is, indeed, a filter. A limit with respect to $\F_F$ coincides
with ordinary limit.
\begin{definition}
A filter $\F$ is called to be ultrafilter iff for any set $A\subseteq \Naturales$ one has
$A\in \F$ or $\Naturales\backslash A\in \F$.
\end{definition}
\begin{theorem}
A bounded sequences has a limit with respect to an ultrafilter. This limit is unique.
\end{theorem}
{\bf Example} For $i\in\Naturales$ let  $\F_i=\{A\subseteq \Naturales\ | \ i\in A\}$. It is an ultrafilter.
Such an ultrafilter is called proper for $i$.
One can check that
$\lim_{\F_i} a_n=a_i$. So, limits with respect to a proper ultrafilter are not interesting.
\begin{proposition}
An ultrafilter $\F$ is proper (for some $i\in \Naturales$) if and only if it contains a finite set.
\end{proposition}
This proposition implies that an ultrafilter is non-proper if and only if it is an extension
of the Frech\'et filter $\F_F$.
On the other hand, it follows  from the Zorn lemma that any filter can be extended to an ultrafilter.
\begin{proposition}
There is an ultrafilter $\F\supset \F_F$. Any such an ultrafilter is non-proper.
\end{proposition}
\subsection{Appendix II}

The Marriage Lemma of P. Hall
(see, for instance, \cite{Ryser}) is formulated as follows.
\begin{lemma} \nonumber \label{mariage}
For an indexed collections of finite sets $F_1, F_2,\ldots, F_k$ the following conditions
are equivalent:
\begin{itemize}
\item there exists an injective function $\alpha : \{1,2,...,k\}\to\bigcup\limits_{i=1}^k F_i$
such that $\alpha(i)\in F_i$;
\item For all $S\subseteq \{1,2,\ldots,k\}$ one has $|\bigcup\limits_{i\in S} F_i|\geq |S|$.
\end{itemize}
\end{lemma}

address: {\small IICO-UASLP, Karakorum 1470, Lomas 4a 78210, San
Luis Potosi, S.L.P. MEXICO

email: valentin@cactus.iico.uaslp.mx, glebsky@cactus.iico.uaslp.mx


\begin{thebibliography}{xx99}

\bibitem{AG}
V. Afraimovich and L. Glebsky, Measures of $\epsilon$-complexity,
{\it Taiwanese J. of Math.} {\bf 9} (2005), 397-409

\bibitem{AG1}
V. Afraimovich and L. Glebsky, Complexity, fractal dimensions and topological entropy
in dynamical systems,
in:
Chaotic Dynamics and Transport in Classical and Quantum Systems,
P. Collet et al. (Eds), Kluwer Academic Publishers (2005), 35-72.

\bibitem{AZ}
V. Afraimovich and G.M. Zaslavsky, Space-time complexity in
Hamiltonian dynamics, {\it Chaos} {\bf 13} (2003), 519,532.


\bibitem{BHM}
F. Blanchard, B. Host, and A. Maass, Topological complexity,
{\it Ergod. Theory Dyn. Syst.} {\bf 20} (2000), 641-662.


\bibitem{B}
M. Boshernitzan, A condition for minimal interval exchange maps to be
uniquely ergodic, {\it Duke Math. J.,}
{\bf 52} (1985), 723-752.

\bibitem{bu} N. Bourbaki, {\it Elements of mathematics. General topology. Part 1.} Hermann, Paris, 1966.

\bibitem{Bo}
R. Bowen, Topological entropy for noncompact sets, {\it Trans. AMS}
{\bf 84} (1973), 125-136.

\bibitem{Bo1} R. Bowen, Entropy for endomorphisms and homogeneous space, {\it Tran. AMS}
{\bf 153} (1971), 404-414.


\bibitem{Fe}
S. Ferenczi, Complexity of sequences and dynamical systems, {\it Discrete Math.}
{\bf 206} (1999), 145-154.


\bibitem{GG} L.Yu. Glebsky, E.I. Gordon and C.J. Rubio, On approximation of unimodular
groups by finite quasigroups, {\it Illinois Journal of Mathematics}
{\bf 49} (2005), 17-31

\bibitem{Halm} P.R. Halmos, {\it Measure theory}, Springer-Verlag, New York, 1974, MR 0033869

\bibitem{KT}
A.N. Kolmogorov and V.M. Tikhomirov, $\epsilon$-entropy and
$\epsilon$ capacity of sets in functional spaces, {\it Usp. Mat. Nauk}
{\bf 14} (1959), 3-86.


\bibitem{Ryser} H. J. Ryser, Combinatorial mathematics,
{\it The Carus Mathematical Monographs}, {\bf 15}
The Mathematical Association of America, 1963.

\bibitem{T} F.Takens, Distinguishing deterministic and random systems, in:
Nonlinear Dynamics and Turbulence, G.I.Barenblatt, G.Iooss, D.D.Joseph, eds.,Pitman (1983) 314-333.

\bibitem{V}
W.A. Veech, The metric theory of interval exchange transformations III,
The Sah Arnoux Fathi invariant, {\it Amer. J. Math}
{\bf 106} (1984), 1389-1422.

\bibitem{ZA}
G.M. Zaslavsky and V. Afraimovich, Working with complexity functions in:
Chaotic Dynamics and Transport in Classical and Quantum Systems,
P. Collet et al. (Eds), Kluwer Academic Publishers (2005), 78-85.


\end{thebibliography}
\end{document}